\documentclass{dcds}
\usepackage{relsize}

\usepackage{amsmath,amssymb}
  \usepackage{paralist}
 \usepackage{psfrag, subfigure}
  \usepackage{graphics} 
  \usepackage{epsfig} 
\usepackage{graphicx}  \usepackage{epstopdf}
 \usepackage[colorlinks=true]{hyperref}
\hypersetup{urlcolor=blue, citecolor=red}

\def\currentvolume{}
 \def\currentissue{}
  \def\currentyear{}
   \def\currentmonth{}
    \def\ppages{X--XX}
     \def\DOI{}

\newtheorem{thm}{Theorem}[section]
\newtheorem{clly}{Corollary}[section]

\newtheorem{lema}{Lemma}[section]
\newtheorem{prop}{Proposition}[section]

\theoremstyle{definition}
\newtheorem{defi}{Definition}[section]
\newtheorem{rk}{Remark}[section]

\newcommand{\ep}{\varepsilon}

\usepackage{tikz-cd}

\newcommand{\pf}{{\flushleft{\bf Proof: }}}
\newcommand{\N}{\mbox{$\mathbb{N}$}}
\newcommand{\T}{\mbox{$\mathbb{T}$}}
\newcommand{\Z}{\mbox{$Z\!\!\!Z$}}

\newcommand{\R}{\mbox{$I\!\!R$}}

\numberwithin{equation}{section}

\newcommand{\transv}{\mathrel{\text{\tpitchfork}}}
\makeatletter
\newcommand{\tpitchfork}{%
  \vbox{
    \baselineskip\z@skip
    \lineskip-.52ex
    \lineskiplimit\maxdimen
    \m@th
    \ialign{##\crcr\hidewidth\smash{$-$}\hidewidth\crcr$\pitchfork$\crcr}
  }%
}
\makeatother

\begin{document}

\title[ Robust transitivity of singular endomorphisms ]{ A persistently singular map of $\mathbb{T}^n$ that is $C^2$ robustly transitive but is not $C^1$ robustly transitive.}

\author[Juan C. Morelli]{}


%

\email{jmorelli@fing.edu.uy }

\maketitle

\centerline{\scshape  Juan Carlos Morelli Ramírez$^*$}
\medskip
{\footnotesize
 \centerline{Universidad de La Rep\'ublica. Facultad de Ingenieria. IMERL.}
   \centerline{ Julio Herrera y Reissig 565. C.P. 11300.}
   \centerline{ Montevideo, Uruguay.}}

\bigskip

 \centerline{(Communicated by )}

\begin{abstract}

 Let $\mathcal{E}$ be the set of endomorphisms defined on the $n$-torus. We exhibit an example of a map in $\mathcal{E}$ that is robustly transitive if $\mathcal{E}$ is endowed with the $C^2$ topology but is not robustly transitive if $\mathcal{E}$ is endowed with the $C^1$ topology.

\end{abstract}\footnote{This work was partially financed by ANII of Uruguay}

\section{Introduction}

Whenever we think about dynamical systems' properties almost inevitably come to mind the concepts of {\textit{stability}} and \textit{robustness}. Loosely speaking, we can say that stability implies same dynamics for maps sufficiently close to each other, and robustness implies the same behavior relative to a specifical property for maps sufficiently close to each other. These are both of most importance in the study of any dynamical system.\\

   This work in particular is centered in the study of robust transitivity, meaning by \textit{transitive} the existence of a forward dense orbit of a point. This may seem at first sight as an unexciting topic since a fair amount of results concerning robust transitivity are known. Nonetheless, the aimed class of maps, the singular endomorphisms about which little to nothing is known; as well as taking on the high dimensional context are undoubtedly a fresh approach to the subject.\\

   To set ideas in order we list up the most relevant known results about the topic.\\
   We begin summing up the most studied case: robust transitivity of diffeomorphisms. The image provided by known results is fairly complete. Concerning surfaces, \cite{m2} shows that robust transitivity implies Anosov diffeomorphism and manifold $\T^2$; while in $dim(M)=n$ manifolds, in \cite{bdp} is proved that robust transitivity implies a dominated splitting. \\ Going further there is robust transitivity of regular endomorphisms (not globally but locally invertible). The image we have about these is somewhat less complete: we know that volume expanding is a necessary but not sufficient condition for $C^1$ robust transitivity according to \cite{lp}.\\ Carrying on, at last there is the least studied case, robust transitivity of singular maps (non empty critical set). Until 2013 nothing had ever been written on the topic. It was on that year when \cite{br} showed the first example of a $C^1$ transitive singular map. The second example was given only in 2016 by \cite{ilp}, they show a $C^1$ robustly transitive map with a persistent critical set. Nothing more than these two examples was known until that time. \\ Either so, there have been recent further advances on the topic: in 2019 \cite{lr1} and \cite{lr2} set the \textit{state of the art} proving that partial hyperbolicity is a necesary condition for robust transitivity of singular surface endomorphisms, that the only surfaces that support them are $\T^2$ and the Klein bottle, and that they belong to the homotopy class of a linear map with an eigenvalue of modulus larger than one. \\

   Now, about the present work, the construction carried on in \cite{ilp} allows the existence of an endomorphism of $\T^2$ with persistent critical set which is $C^2$ but not $C^1$ robustly transitive, a result appearing in \cite{ip} which this article generalizes to the higher dimensional torus $\T^n$. It's worth to mention that the proof is inspired in the preceding ideas but with a significantly simpler approach and construction.

    \subsection{Sketch of the Construction.} Starting from a matrix with integer coefficients of absolute value larger than one we build and endomorphism of $\T^n$ that presents a persistent critical set and admits a field of unstable cones. We choose a critical point and a neighborhood of it and perform a perturbation there. It provides with a new map that collapses an open set to an invariant hyperplane so it can not be $C^1$ transitive. Next, we prove that the map is $C^2$ robustly transitive by following curves whose velocities lie inside the unstable cones until they escape the critical region and then apply the classical argument for robust transitivity using open sets since no open set collapses in the $C^2$ topology. \\ The reader is also provided with a thorough description of the critical set and the critical points.

\section{Preliminaries}

We begin recalling some basic definitions. We assume the reader to be familiar with the concepts of real manifold and submanifold, atlas, chart, tangent vector and tangent space, differentiable map and differential of a map, etc. For more details about the contents of this section the reader might refer itself to \cite{gg} or \cite{kh}.

\subsection{Dynamical}
Let $f:M \rightarrow M$ a differentiable endomorphism. The \textbf{orbit} of $x \in M$ is $\mathcal{O}(x)=\{ f^n(x) , n \in \N \}$ and $f$ is \textbf{transitive} if there exists a point $x \in M$ such that $ \overline{{\mathcal{O}(x)}}=M$.

\begin{prop}\label{equi}
  If $f$ is continuous then are equivalent:
  \begin{enumerate}
    \item  $f$ is transitive.
    \item  For all $U, V$ open sets in $M$, exists $n \in \N$ such that $ f^n(U) \cap V \neq \emptyset$.
    \item  There exists a residual set $R$ (countable intersection of open and dense sets) such that for all points $  x \in R: \overline{\mathcal{O}(x)}=M$.
  \end{enumerate}
\end{prop}

\begin{defi}
   $f$ is $C^k$-\textbf{robustly transitive} if there exists $\ep >0$ and a neighborhood $\mathcal{U}_{(f,\ep)}$ of $f$ in the $C^k$ topology such that $g$ is transitive for all $ g \in \mathcal{U}_f$.
\end{defi}

\subsection{Geometrical}

For the rest of the preliminaries $M$ and $N$ will denote real manifolds, compact, connected and without boundaries such that $dim(M) \geq  dim(N)$, $x \in M, y \in N$ and a differentiable map $f:M \rightarrow N$.\\
  We say that $x $ is a \textbf{regular point} for $f$ if the differential at $x$, $D_x f$ is surjective. Or equivalently, if the rank of the Jacobian matrix of $f$ at $x$ satisfies $rk(D_x f) = dim(N)$. We say that $y $ is a \textbf{regular value} of $f$ if $\forall x \in f^{-1}(y)$, $x$ is a regular point and we say that $x \in M$ is a \textbf{critical point} or \textbf{singularity} for $f$ if $D_x f$ is not surjective. Equivalently, if $rk(D_x f) < dim(N)$. The \textbf{critical set} of $f$ is $ S_f = \{ x \in M / rk(D_x f) < dim(N) \} $ and $y $ is a  \textbf{critical value} if it is not a regular value. 
  \begin{rk}
  If $dim(M)=dim(N)$ then the definitions of regular point and critical point are equivalent to the determinant of the Jacobian $det(D_xf) \neq 0$ or $det(D_xf) = 0$ respectively.
\end{rk}

 \subsubsection{Singularities}

 We continue with a brief overview of singularity theory. Surprisingly, singularities present a nice geometrical behavior; they can be grouped under 'types' of singularities and these types can as well be grouped as submanifolds of $M$.

\begin{defi}
  Let $f_1:M\rightarrow N$ and $f_2:M\rightarrow N$ maps such that $x_1 \in S_{f_1}$ and $x_2 \in S_{f_2}$.
  We say $x_1$ and $x_2$ are \textbf{same type singularities} and denote it by $(f_1,x_1) \sim (f_2,x_2)$ if there exist neighborhoods $U_i$ of $x_i$, $V_i$ of $f_i(x_i)$, $i \in \{1,2\}$ and two diffeomorphisms $h_1:U_1 \rightarrow U_2$ and $h_2: V_1 \rightarrow V_2$ such that the following diagram commutes:
 \begin{Large}{ $$\begin{tikzcd}
U_1 \arrow{r}{h_1} \arrow[swap]{d}{{f_1}_{|U_1}} & U_2 \arrow{d}{{f_2}_{|U_2}} \\
V_1 \arrow{r}{h_2} & V_2
\end{tikzcd}$$}
\end{Large}
\end{defi}

Clearly, $\sim$ is an equivalence relation. \\
The description of all singularities belonging to the same equivalence class receives the name of \textbf{normal form} of the singularity, and for simplicity they are classified for some map $g:\R^m \rightarrow \R^n$ with a singularity at $x=0$. Then, we say that $f:M \rightarrow N$ has \emph{that} type of singularity at $x \in M$ if $(f,x) \sim (g,0)$. \\
This description carries on with the disadvantage of allowing equivalence classes being too many, even infinite. Considering this, equivalence classes are grouped once again under a criterion called \textbf{Thom-Boardman}. Every differentiable map can be approximated by another with a finite number of Thom-Boardman singularities, and for this description it holds that every group of classes is, as a set, a submanifold of $M$. \\
In 1955, Thom proposed in \cite{t} the afore mentioned criterion for the singularities of a map $f \in C^\infty(M,N)$ defining $S_{k_1}(f)=\{ x \in M / dim (ker(D_xf))={k_1}\}$. Assuming this set to be a submanifold of $M$, then $S_{k_1,k_2}(f)=S_{k_2} (f_{|S_{k_1}})$ can be defined and so on. The non-increasing r-tuple $(k_1,k_2,...,k_r)$ is called to be the \emph{symbol} of the singularity and it characterizes it.\\ Unfortunately Thom could not prove the submanifold structure, leaving a blank in the theory until Boardman solved the problem in 1967 in his paper \cite{b}. He took jets of differentiable maps (loosely speaking Taylor expansions) which allow to define the $S_k$ as submanifolds of jets spaces, bringing along as well the local expressions of the action of the map in a neighborhood of the singularity (which happen to be the explicit analytical expressions of their normal forms). This sealed the acceptance of the classification (another very important consequence of Boardman works on singularities is the proof that stable maps are not dense in the set of smooth maps).

To end with preliminaries concerning singularities we will only mention that Mather (1971) and Morin (1972) gave different approaches to the Thom-Boardman classification by means of algebraic geometry and complex analysis, and that even extremely useful, the Thom-Boardman classification is not free of pathological behavior as Porteous (1972) showed by exhibiting a map $f:\R^5 \rightarrow \R^5$ with $S_2$ singularities in the closure of $S_{1,1,1,1}$ with very different qualitative behaviors. The boundaries of the Thom-Boardman classification lie in the fact that the closure of the union of the singular submanifolds need not be a manifold itself.

\section{A singular endomorphism $f$ of $\T^n$.}

We proceed now to the construction of the main map of this work which we call $f$. From $f$ we will be able to  define a new map $H$ that has the properties announced in the title of the article.

\subsection{Construction of $f$.}

 Consider the $n$ dimensional torus $\T^n$ endowed with the standard riemannian metric and let $\widehat{A} \in \mathcal{M}_n(\Z)$ be the diagonal matrix suggested below.
 $$\widehat{A}= \left(\begin{array}{ccccc}
 8 & 0 & 0 & \cdots & 0\\
0 & 2 & 0 & \cdots & 0 \\
\vdots & \vdots &2 & \cdots & 0 \\
\vdots & \vdots & \vdots & \ddots & \vdots \\
0 & \cdots & 0 &0 &2\\
\end{array}%
\right). $$
Notice the construction could be carried on with any pair of integers $\lambda$ and $\mu$ such that $|\lambda|>|\mu|>1$. The choice of $8$ and $2$ is made in the sake of simplicity and for a better understanding of the contents to follow.\\ $\widehat{A}$ defines a regular endomorphism 
$ A:\T^n\rightarrow\T^n /A(x_1,...,x_n)= (8x_1,2x_2,...,2x_n) $. Observe that  $p=(\frac{1}{4},0,...,0,\frac{1}{4})\in\T^{n}$, that $Ap =\left( 0, 0, ..., 0, \frac{1}{2} \right) $ and that $ A^2p=0$. \\ This point $p$ is the center of a ball where a perturbation will be performed in order to obtain the map $f$ we seek. To construct the perturbation we need to fix a series of technical parameters that will define $f$. The choice to do it at the beginning and all of them at the same time is in expectance of avoiding darkness in the construction and of that it will be clear how they depend on each other.\\

Start with $r>0$ satisfying the following conditions:

\begin{itemize}\label{p1}
 \item $A(B_{(p,r)})\cap \overline{B_{(p,r)}}=\emptyset $,
 \item $ A^{-1}(B_{(p,r)})\cap \{ x\in \R ^{n} / x_{n}=0 \}=\emptyset$,
 \item $B_{(p,r)}\cap B_{(A.p ,r)}=B_{(A.p ,r)}\cap B_{(0,r)}= \emptyset.$ 
 \end{itemize}

\begin{figure}[ht]
\begin{center}
\subfigure[]{\includegraphics[scale=0.33]{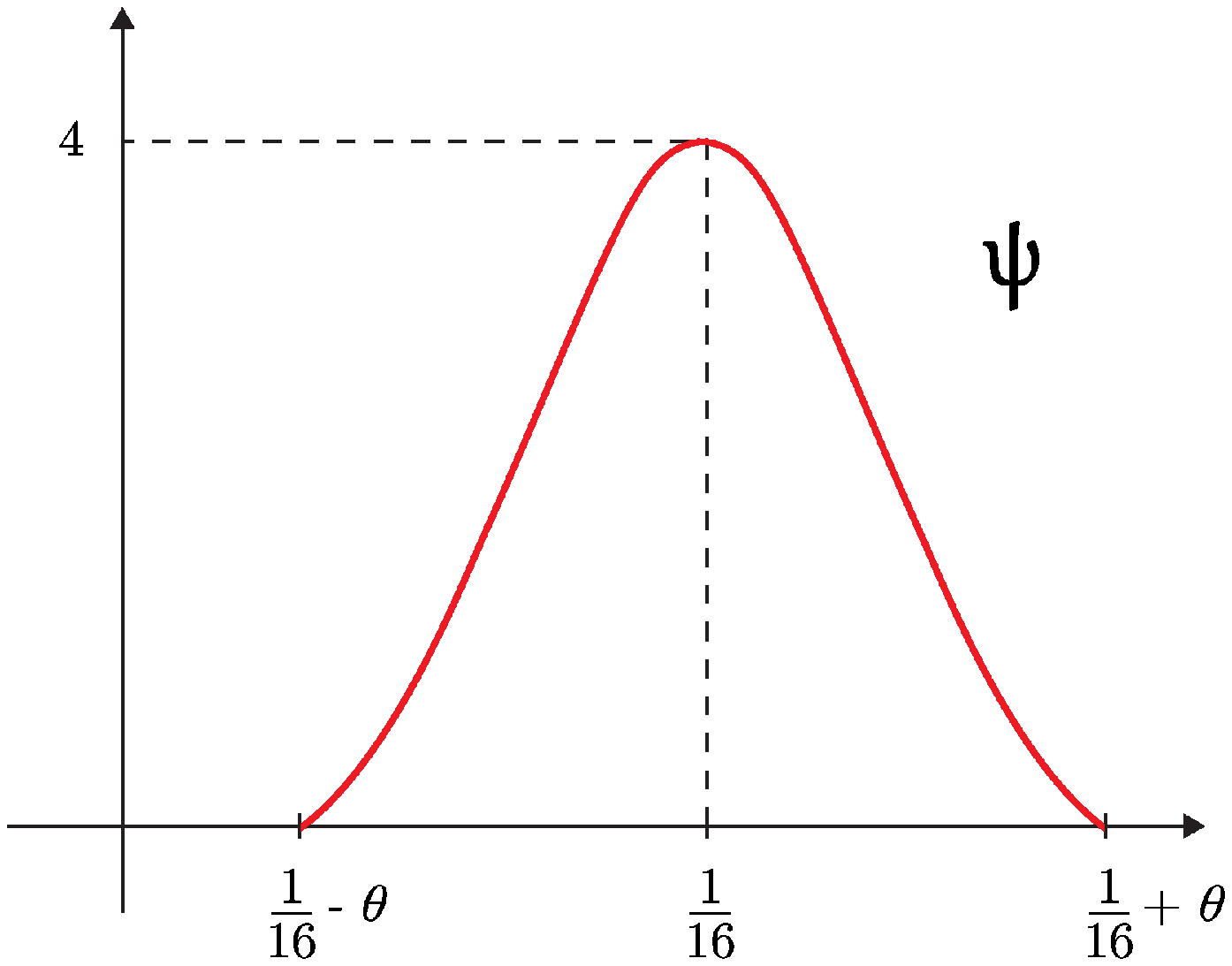}}
\subfigure[]{\includegraphics[scale=0.37]{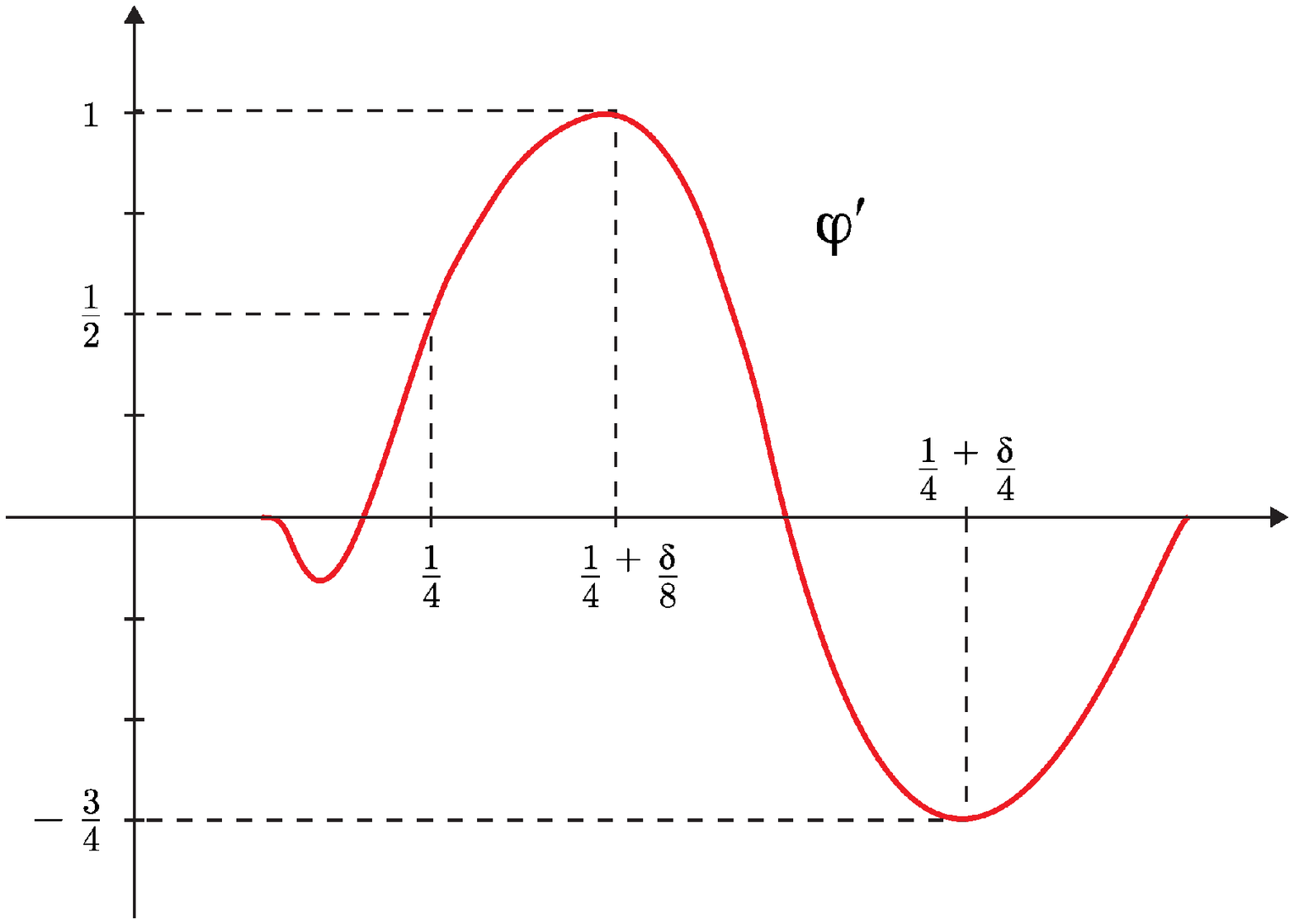}}
\caption{Graphs of $\psi$ and $\varphi'$}\label{figura11}
\end{center}
\end{figure}

We choose $r$ like this so as there will be no points in $B_{(p,r)}$ that remain inside it under forward iteration by $A$ or coming from the hyperplane $\{x_n=0\}$ under iteration by $A$.\\
Now that $r$ is set, fix a second parameter $\theta$ such that  $0< \theta  < \frac{r}{2}$  and define a function $\psi :\R\to\R$ of class $C^{\infty}$ with an only critical point at $\frac{1}{16}$, with $\psi (\frac{1}{16})=4$ and $\psi ^{'}(\frac{1}{16})=\psi ^{''}(\frac{1}{16})=0$, with $\psi (x)=0$ for all $ x$ in the complement of $(\frac{1}{16}-\theta ,\frac{1}{16}+\theta  )$;  with an axis of symmetry in the line $x=\frac{1}{16}$ as shown in Figure \ref{figura11} (a) .

Choose now a third arbitrary real parameter $a \in (0,\frac{3}{7})$ (which will play a key role shortly on in Section \ref{dinamf}) and finish setting a last parameter $\delta $, with  $0<\delta < 2\theta $ verifying the following condition: since the derivative of $\psi$ is bounded once $\theta$ has been fixed, call $ M := M(\theta) = max_{x \in \R}  \{ | \psi'(x) | \} $ and impose on $\delta$ that $2.M.r.\delta.(1+a)<a.$  \\

Now that all the parameters have been fixed, consider a smooth (class $C^{\infty}$) function $\varphi:\R\to\R$ such that:

\begin{itemize}
\item $\varphi' $ is as in Figure \ref{figura11} (b),
\item $\varphi (\frac{1}{4})=0$, $\varphi' (\frac{1}{4})=\frac{1}{2}$,  $\varphi^{''}(\frac{1}{4})\neq 0$ , $\varphi' (\frac{1}{4} +\frac{\delta}{8}  )=1$,
\item $\varphi (x)=0$ for all
$ x \notin [\frac{1}{4}-\frac{\delta}{4},\frac{1}{4}+\frac{3\delta}{4} ]$.
\item $-\frac{3}{4} \leq \varphi'(x) \leq 1$ for all $x \in [\frac{1}{4}-\frac{\delta}{4},\frac{1}{4}+\frac{3\delta}{4} ]$.
\end{itemize}
\begin{rk}
  $max \{ |\varphi (x)|: \ x \in \R \}\leq \delta.$
\end{rk}

We are now in condition to define the perturbation of $A$ that depends on $r$,$\theta$, $a$ and $\delta$ which by simplicity we call only $f$ and is  \begin{equation}\label{mapaefe}
           f:= f_{r,\theta,a, \delta}:\T^n\to \T^n / f(x_1,...,x_n)= \left( 8x_1,2x_2,..., 2x_n-\varphi(x_n).\psi \left( \sum_{k=1} ^{n-1} x_k^2 \right) \right).
         \end{equation}

To make the reading easier we will denote it as $ f(x_1,...,x_n)= \left( 8x_1,2x_2,..., 2x_n-\varphi.\psi  \right)$ omitting the evaluations appearing on the definition.

\begin{rk}\label{obs1}. It is straightforward that:\\
  1- $f \in C^{\infty}$.\\
  2- $f_{|B^c_{(p,r)}}=A$. 	\\
  3- $ f(p)= Ap =\left( 0, 0, ..., 0, \frac{1}{2} \right)$  and  $f^2(p)= A^2p=0$.\\
  4- The differential $D_x f$ at $x=(x_1,...,x_n)$ is

\begin{equation}\label{ecuacion1}
D_xf =
\left( \begin{array}{ccccc}
8 & 0 & \cdots & 0 & 0 \\
0 & 2 & \cdots & 0 & 0 \\
\vdots & \vdots & & \vdots & \vdots \\
0 & 0 & \cdots & 2 & 0 \\
-2.x_1.\varphi.\psi' & -2.x_2.\varphi.\psi'& \cdots & -2.x_{n-1}.\varphi.\psi' & 2-\varphi'.\psi \\
\end{array}
\right).
\end{equation}\\
  \end{rk}

\subsection{Dynamics of $f$.}\label{dinamf}
In this subsection we prove that the dynamical behavior of $f$ is given by the existence of strong unstable cones at every point in the direction of the canonical first coordinate. We recall some definitions from \cite{kh} first.

\begin{defi}\label{cono}
  Let $x \in M$, we call  \textbf{cone} of parameter $a$ and vertex $x$  to  $$C^{u}_{a} (x)=\left\{ (v_{1},...,v_{n}) \in T_x M / \frac{\Vert (v_{k+1},...,v_{n}) \Vert}{\Vert (v_{1}, v_2,...,v_k) \Vert}<a \right\}$$ for each $k \in \N \cap [1,n-1]$.
\end{defi} In this case we say the unstable cone is of \textbf{index} $n-k$.
\begin{defi}
  We say a map $f$ admits an \textbf{unstable cone} of parameter $a$ and vertex $x \in M$ if exists $C^u_a(x) \subset T_xM$ such that $\overline{D_x f(C^u_a(x)) }\setminus \{0\} \subset C^u_a(f(x))$.\\
  We say $f$ admits unstable cones of parameter $a$ if it admits an unstable cone of parameter $a$ and vertex $x$ at every point $x \in M$.
\end{defi}
\begin{prop}
{\textsl{Existence of unstable cones is $C^1$ robust}}:\\
   If $f$ admits unstable cones of parameter $a$ then there exists a neighborhood $\mathcal{U}_f \in C^1$ of $f$ such that $\forall g \in \mathcal{U}_f: g$ admits unstable cones of parameter $b$ with $b \leq a$.
\end{prop}

In the context of diffeomorphisms, the existence of unstable cones is equivalent to a weak form of hyperbolicity called dominatted splitting. The subspace generated by the vectors in the denominator of the cone is regarded as the unstable space. Though this condition is not an equivalence for endomorphisms, the dynamical behavior of a small parameter cone is virtually the same as that of a strong unstable direction.\\ As was stated above, the most relevant dynamical feature that our map $f$ has is the existence of strong unstable cones at every point. We give the proof of this assertion below. Observe that the choice of the technical third parameter $a$ before fixing $\delta$ in the construction of $f$ is what we need for the cones to exist. Recall that we imposed on $\delta$ the condition $2.M.r.\delta.(1+a)<a$ for the parameter $a$ and the norm $M$ of the derivative of $\psi$. Observe also that for every $\delta_0$ such that $0 < \delta_0 < \delta$ the claims stated up next for $\delta$ are satisfied.\\

We will make use of the following notations in the proof of the lemma ahead and also in the following sections:

\begin{itemize}
\item If $v=(v_{1},v_{2},...,v_{n})$ and $h<n$, then $\tilde{v}_{h}=(v_{1},v_{2},...,v_{h})$ and $\tilde{v}_{n-1} = \tilde{v}$.
\item If $v=(v_{1},v_{2},...,v_{n})$ and $w=(w_{1},w_{2},...,w_{n})$ then $ \langle v,w\rangle_{h} = \sum_{k=1}^{h} v_{k}.w_{k} $.
\end{itemize}

\begin{lema} \label{3}
\textbf{Existence of unstable cones for $f$.} \end{lema}

Claim: Given $a \in  (0, \frac{3}{7}) $:
\begin{enumerate}
\item $ C^{u}_{a} (x)=\lbrace (v_{1},...,v_{n})/ \frac{\Vert v_{2},...,v_{n} \Vert}{\vert v_{1}\vert}<a \rbrace$ satisfies $\overline{D_{x}f(C^{u}_{a}(x))} \setminus \{0\} \subseteq C^{u}_{a}(f(x))$.\\
\item For all $ v \in C^{u}_{a} (x)$ holds that $ \Vert D_{x}f (v) \Vert > 7 \Vert v \Vert $.
\end{enumerate}

\pf
\begin{enumerate}
\item  By Equation \ref{ecuacion1} we have $\forall v=(v_{1},v_{2},...,v_{n}) \in C^{u}_{a} (x):$

$$u (u_1,..,u_n):=D_{x}f (v) = (8.v_{1},2.v_{2},...,-2.\langle x,v\rangle_{n-1}.\varphi.\psi' + v_{n}.(2-\varphi'.\psi)). $$

Performing calculations we have  $$ \frac{\Vert u_{2},...,u_{n} \Vert}{\vert u_{1}\vert}  = \frac{\Vert (2.v_{2},...,2.v_{n-1}, -2.\langle x,v\rangle_{n-1}.\varphi .\psi'+ v_{n}.(2-\varphi' .\psi) \Vert}{\vert 8.v_{1} \vert} \leq$$ $$ \leq \frac{\Vert 2.v_{2},...,2.v_{n-1} \Vert}{\vert 8.v_{1}\vert} + \frac{2.\Vert \tilde{x} \Vert.\Vert \tilde{v} \Vert.\vert \varphi \vert.\vert \psi' \vert}{\vert 8.v_{1} \vert} + \frac{\vert 2-\varphi' .\psi \vert. \vert v_{n}\vert}{\vert 8.v_{1} \vert} < \frac{a}{4} + 2.r.(\frac{1+a}{8}).\delta.M +\frac{5.a}{8} < a; $$ \\ where in the first inequality we use triangular and Cauchy-Schwarz; and in the second one we use:
\begin{itemize}
\item $ v \in C^{u}_{a} (x)$,
\item $\Vert \tilde{x} \Vert \leq \Vert x \Vert < r$,
\item $ \frac{\Vert \tilde{v} \Vert}{|8.v_1|}\leq \frac{|v_1|+\Vert v_2,...,v_{n-1} \Vert}{|8.v_1|}\leq \frac{1+a}{8}$,
\item $ \vert \varphi \vert < \delta $,
\item $ \vert \psi' \vert < M $,
\item $ \vert 2-\varphi' .\psi \vert \leq 5 $ since $ \frac{-3}{4} \leq \varphi' \leq 1 $ and $ 0 \leq \psi \leq 4$.

\end{itemize}
And in the third one we use the condition imposed over $\delta$.$\Box$ \\
\item  $\forall v=(v_{1},v_{2},...,v_{n}) \in C^{u}_{a} (x)$: $$\left( \frac{\Vert D_{x}f (v) \Vert}{7.\Vert v \Vert}\right) ^{2} = \frac{64.v_{1}^{2}+4.\sum_{j=2}^{n-1}.v_{j}^{2}+(-2.\langle x,v\rangle_{n-1}.\varphi .\psi'+v_{n}.(2-\varphi' .\psi))^{2}}{49.\sum_{j=1}^{n}v_{j}^{2}}  \geq$$
$$ \geq  \frac{64.v_{1}^{2}}{49.\sum_{j=1}^{n}v_{j}^{2}} > \frac{64 }{49(1+a^{2})}> \frac{58 }{49(1+a^{2})} > 1, \forall a \in \left(0,\frac{3}{7}\right).\Box $$
\end{enumerate}

\begin{clly}\label{coroconos}
 There exists a neighborhood $\mathcal{U}_f \in C^1$ of $f$ such that for all $g \in \mathcal{U}_f:g$ admits unstable cones satisfying for all $ v \in C^{u} (x)$ : $ \Vert D_{x}g (v) \Vert \geq 7 \Vert v \Vert $.
\end{clly}

\begin{rk}
  $\vec{e_1}$ can be regarded as a strong unstable direction for $f$.
\end{rk}

\subsection{The critical set of $f$}
Recall that $S_f=\{x \in \T^{n} /det ( D_xf )=0\}$ is the critical set of $f$. Equation \ref{ecuacion1} provides $det(D_xf)=8 \cdot (2)^{n-2} \cdot (2-\varphi'.\psi) $ which translates into $$S_f=\{x \in \T^{n} / 2-\varphi'.\psi =0\rbrace.$$

\begin{rk}\label{rkcritset}
  1-  The point $p \in S_{f} $, so $S_f\neq \emptyset$.\\
2- The last column in the Jacobian (\ref{ecuacion1}) is null at every point belonging to  $S_f $. Hence, for every point $ x \in S_f$ the kernell $ker(D_x f)$ is generated by the last canonical vector $\{\vec{e_n}\} $.\\
3- \textbf{The critical set $S_f$ is persistent} in the sense that for all $k \geq 1$ there exists a $C^k$ neighborhood of  $f$ for which all of its maps have nonempty critical set. Take the points $q_1=(\frac{1}{4},0,...,0, \frac{1}{4} + \frac{\delta}{4})$ and $q_2=(\frac{1}{4},0,...,0, \frac{1}{4} + \frac{\delta}{8})$ both in $B_{(p,r)}$. Evaluate determinants to obtain \textit{det}$(D_{q_1}f)= 5. (8^{n-1})$ and  \textit{det}$(D_{q_2}f)= -2.( 8^{n-1})$. Therefore it exists $\mathcal{U}_f \in C^{1}$ such that $\forall g \in \mathcal{U}_f$: $S_g\neq\emptyset$.
\end{rk}

 \subsubsection{Classification of $S_f$}\label{Sfessubvariedad}

The objective of this subsection is to prove that the critical set $S_f$ is a submanifold of $\T^n$ of codimension 1. We begin recalling some basic definitions and theorems from differential geometry that will be necessary to carry on with our proof. Again, for more insight on these concepts the reader might refer itself to \cite{gg}. Remember that $M$ and $N$ denote differentiable manifolds with $dim(M) \geq dim(N)$, $x \in M, y \in N, W \subset N$ a submanifold and $f:M \rightarrow N$ a differentiable map.

\begin{defi}\label{transve}
We say that $f$ is \textbf{transversal} to $W$ at $f(x) \in N$ , and we denote it by $f \transv W$ at $f(x)$ if any of the two following conditions hold: $ f(x) \notin W$ or $T_{f(x)}W\oplus D_xf(T_xM)=T_{f(x)} N $.
\end{defi}

\begin{defi}
  We say that $f \transv W$ if for all $x \in f^{-1}(W)$ holds $f \transv W$ at $f(x)$.
\end{defi}

\begin{prop}
  Transversality is a \textbf{$C^1$-open} property ($C^1$-stable) i.e. $f \transv W$ then there exists a neighborhood  $\mathcal{U}_f \in C^1$ of $f$ such that $\forall g \in \mathcal{U}_f: g \transv W$.
\end{prop}

\begin{thm}\label{subvartransversal}
{\bf{Preimage of transversal submanifold}}:
Let $f:M\rightarrow N$ be a differentiable map. $W$ a submanifold of $N$ and $f \transv W$.
Then $f^{-1}(W)$ is a submanifold of $M$ of the same codimension that $W$ in $N$.
\end{thm}

\begin{thm}\label{valorregular}
{\bf{Preimage of a regular value}}:
Let $f:M\rightarrow N$ be differentiable and $y \in N$ a regular value.
Then $f^{-1}(y)$ is a submanifold of $M$ with codimension equal to the dimension of $N$. Also, $\forall x \in f^{-1}(y)$ holds $T_x\left( f^{-1}(y)\right) =ker(D_xf)$.
\end{thm}

\vspace{0.2cm}
We start developing now the ideas to prove that $S_f$ is a submanifold of $\T^n$.

 Let $ \mathcal{M}_{n}(\R)$ be the set of matrixes with real coefficients and size $n \times n$, and let $\mathcal{R}_{k}$ be the subset of matrixes of rank $k$. Then $\mathcal{R}_{k}$ is a submanifold of $\mathcal{M}_{n}(\R)$ (\cite{gg}, \textit{chapter 2, proposition 5.3}). In particular, if $k=n-1$ then $\mathcal{R}_{n-1}$ is of dimension $n^2-1$. Define now a map $h$ depending on $f$ by
 \begin{equation}\label{mapah}
   h:= h_{f}:\T^n \rightarrow \mathcal{M}_{n }(\R) / h\left( x\right) = D_{x}f
 \end{equation}

 which assigns to every point $x$ of $\T^{n}$ the Jacobian matrix of $f$ at $x$. Since $f$ is smooth so is $h$. Observe that $h(S_f)\subset \mathcal{R}_{n-1}$ since for all $ x\in S_f$ the rank $rk(D_{x}f)=n-1 $, but also since $h^{-1}\left(\mathcal{R}_{n-1}\right)\subset S_f$ then $S_f=h^{-1} \left(\mathcal{R}_{n-1}\right)$.\\
 Proving that $ h \transv \mathcal{R}_{n-1}$ would imply our claim via Theorem \ref{subvartransversal}. \\Observe that if $x \notin S_f$, then $h(x) \notin \mathcal{R}_{n-1}$ hence $ h \transv \mathcal{R}_{n-1}$. \\
  It is only left to see that the transversality condition holds when $x \in S_f$.\\ By definition $h \transv \mathcal{R}_{n-1}$ at $x \Leftrightarrow T_{h(x)}\left(\mathcal{M}_{n}(\R)\right)=T_{h(x)}\left(\mathcal{R}_{n-1}\right)\oplus \left( D_{x}h \right)(T_{x}(\T^n))$. Or equivalently, $$h \transv \mathcal{R}_{n-1} \mbox{ at } x \Leftrightarrow T_{h(x)}\left(\R^{n^2}\right)=T_{h(x)}\left(\mathcal{R}_{n-1}\right)\oplus \left( D_{x}h\right)(\R^n).$$

     We determine $ T_{h(x)}\left(\mathcal{R}_{n-1}\right)$ and $Im \left( D_{x}h\right)$ when $x \in S_f$ to check that the transversality condition holds.

\begin{lema} \label{lema 1} For every critical point $x$ it holds that $$T_{h(x)}\left(\mathcal{R}_{n-1}\right)=\lbrace (x_1,...,x_{n^2})\in \R^{n^2} /x_{n^2}=0\rbrace. $$ \end{lema} \pf

Consider the determinant map $\xi:\mathcal{R}_{n-1}\rightarrow\R / \xi(M)= det(M)$. Follows immediately that $Im( \xi )=0$. \\ For all $ M \in \mathcal{R}_{n-1}$ the differential of $\xi$ is $D_{M} \xi=(A_{11} , A_{12} ,..., A_{nn} )$ being $A_{ij}$ the adjoint matrix to the element in the $i$-th row and $j$-th column of $M$. Since $M \in \mathcal{R}_{n-1}$, at least one of these adjoints is not null. This implies that the differential of $\xi$ has dimension greater or equal than one for all $ M \in \mathcal{R}_{n-1}$ which means that for all $ M \in \xi^{-1} (0): D_{M} \xi$ is surjective. Consequently, $0$ is a regular value for $\xi$ so according to Theorem \ref{valorregular} the tangent space to $\mathcal{R}_{n-1}$ at $h(x)$ is $T_{h(x)} (\xi^{-1} \lbrace0\rbrace) = ker (D_{h(x)}\xi)$. Now, observe that for every critical point $x$ it is $D_{h(x)}(\xi)=(0,...,0, 8.(2^{n-2}))$ since the restriction $h_{|S_{f}}$ has its last column null. This implies that for every critical point $ x \in S_f : ker (D_{h(x)}\xi) = \lbrace (x_1,x_2,...,x_{n^{2} -1}, 0)\rbrace \subset \R ^{n^2} $. Combine both equalities above to obtain $T_{h(x)}\left(\mathcal{R}_{n-1}\right)=T_{h(x)} (\xi^{-1} \lbrace0\rbrace) = ker (D_{h(x)}\xi)$ which is the thesis. $\Box$

\begin{lema} \label{lema 2} For every critical point $x$ it holds that $$ \{(0,0,...,0,1)\} \subset Im \left( D_{x}h\right) . $$ \end{lema}
\pf
Since $h:\T^n \rightarrow \mathcal{M}_{n \times n}(\R)$ and $f \in C^2$ then $h \in C^1$; therefore, the differential $D_x h \in \mathcal{M}_{n^2 \times n}(\R)$ exists and is continuous. The statement above is equivalent to the jacobian $D_x h$ has a non-zero element in the last row. \\ An explicit calculation provides:

\[
(D_xh)_{ij} =
     \begin{cases}
       \text{$0$} &\quad\text{if $i<n^2-n+1$}\\
       \text{$-2.\varphi.\psi-4.x_{j}^{2}. \psi''$} &\quad\text{if $n^2-n<i<n^2$ and $j=i$} \\
       \text{$-4.x_{i-n^2+n}.x_j.\varphi.\psi''$} &\quad\text{if $n^2-n<i<n^2$ and $j\not\in \lbrace i, n\rbrace$}\\
       \text{$-2.x_{i-n^2+n}.\varphi'.\psi'$} &\quad\text{if $n^2-n<i<n^2$ and $j=n$}\\
\text{$-2x_j.\varphi'.\psi'$} &\quad\text{if $i=n^2$ and $j\neq n$}\\
\text{$-\varphi''.\psi$} &\quad\text{if $i=n^2$ and $j=n$}\\
     \end{cases}
\]
\\
Consider the last row, at value $i=n^{2}$:\\ If $\varphi''\neq 0$, since $\psi >0$ at all points in $S_f$ we find a non-zero element in the last column when $j=n$. \\
If $\varphi''= 0$, looking at the graph of $\varphi'$ (Figure \ref{figura11}) there are only three points where the condition holds. Recall that for all $x \in S_f$ it is $\varphi'.\psi=2$. Discard then both of the points which have negative image under $\varphi'$ since they don't belong to the critical set, so the only one remaining case is at $x_n= \frac{1}{4}+ \frac{\delta}{8}$. \\ It is only left to see then that if $x_n= \frac{1}{4}+ \frac{\delta}{8}$ there exists $j$ such that $-2.x_j.\varphi'.\psi' \neq 0$ or equivalently $x_j.\psi' \neq 0$ since at this value $\varphi'(\frac{1}{4}+ \frac{\delta}{8})=1$. \\ Discuss now two cases considering Equation (\ref{mapaefe}) and Figure \ref{figura11}:\\ If it were $\psi'=0$ then $\sum_{k=1}^{n-1} x_k^2=\frac{1}{16}$. Since $ \psi(\frac{1}{16})=4$ then $\varphi'.\psi = 4$ implies $x \notin S_f$.\\ If it is  $\psi' \neq 0$ and it were for all $j: x_j = 0$, since $ \psi(0)=0$ then $ \varphi'.\psi = 0$ which again implies $x \notin S_f $. \\ Therefore, for all critical points $ x \in S_f$ satisfying $ x_n= \frac{1}{4}+ \frac{\delta}{8}$ there exists $j<n$ such that $ x_j.\psi' \neq 0$.  $\Box$

\begin{thm}\label{teorema 1}Apply Lemmas \ref{lema 1} and \ref{lema 2} along with Theorem \ref{subvartransversal} to obtain Definition \ref{transve} which gives the proof.$\Box$
\end{thm}
 \subsubsection{Geometry of $S_f$ }\label{representacion}

Having established that $S_f$ is a submanifold of $\T^n$ we turn our attention to its geometry. The objective of this section is to prove that the critical set is isomorphic to the product of two spheres $$S_f \simeq S^{n-2} \times S^1.$$ As an introduction, we analyze the problem in low dimensions first and then generalize. Remember that for vectors $v=(v_{1},v_{2},...,v_{n}) \in \R^n$ and $h<n$ we denote $\tilde{v}_{h}=(v_{1},v_{2},...,v_{h})$ and $\tilde{v}_{n-1} = \tilde{v}$.\\

\begin{figure}[ht]
\begin{center}
\includegraphics[scale=0.13]{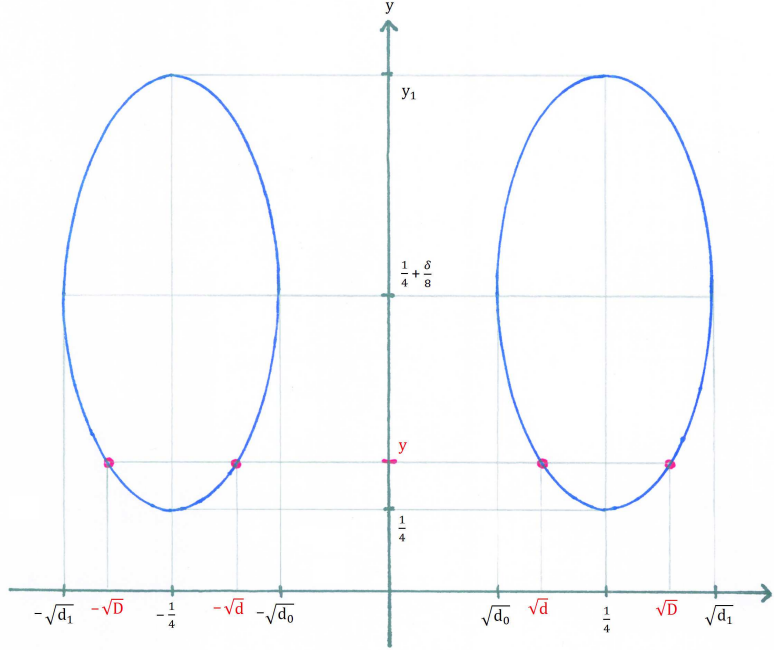}
\caption{$S_f$ for $n=2$.}\label{nigual2}
\end{center}
\end{figure}

 \underline{Case $n=2$}:\\
 This case is fully studied in \cite{ip}. We only observe that the map $f$ takes the form $f:\T^2\to \T^2 / f(x,y)= \left( 8x,2y-\varphi(y).\psi (  x^2 ) \right) $ and give a graph (Figure \ref{nigual2}) showing explicitly that $S_f \simeq S^{0} \times S^1$. \\

\begin{figure}[ht]
\begin{center}
\includegraphics[scale=0.13]{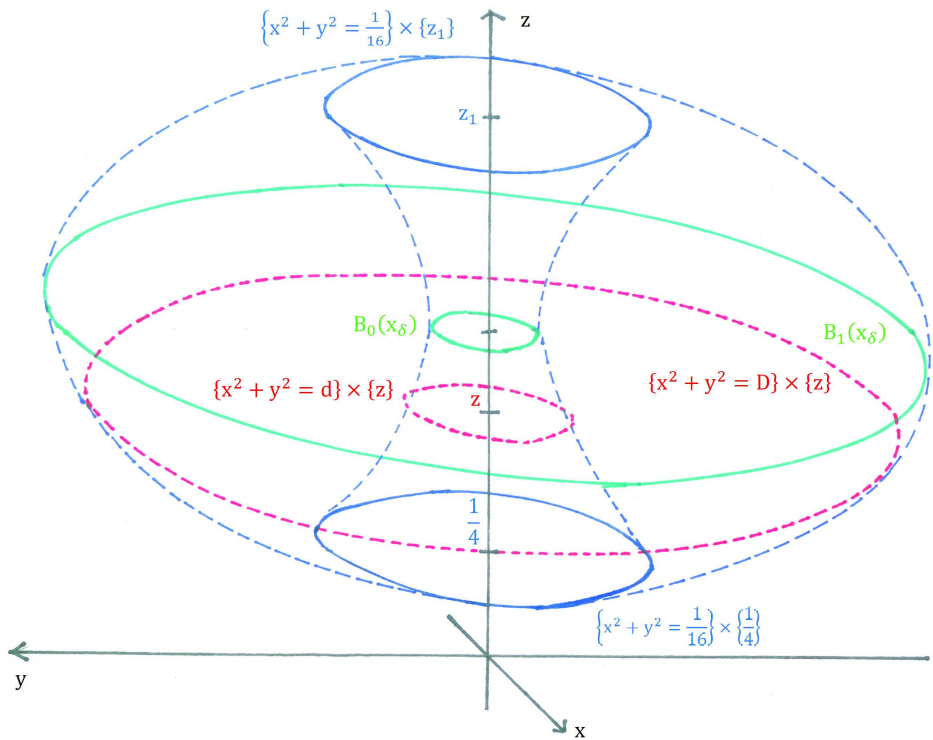}
\caption{$S_f$ for $n=3$.}\label{nigual3}
\end{center}
\end{figure}

 \underline{Case $n=3$}:\\
 Here $f$ is $f:\T^3\to \T^3 / f(x,y,z)= \left( 8x,2y,2z-\varphi(z).\psi (  x^2 + y^2 ) \right)$ and consequently the critical set is $S_f=\{(x,y,z) \in \T^{3} / 2-\varphi'(z).\psi(x^2+y^2) =0\rbrace$. Since $\psi \leq 4$, if $\varphi' < \frac{1}{2}$ then $\psi \cdot \varphi' <2$ which imply for all $(x,y,z) \in S_f : \varphi'(z) \geq \frac{1}{2}$. In turn $\varphi'(z) \in [\frac{1}{2},1]$ and $\psi(x^2+y^2) \in [2,4]$.\\ Go now to the graphs of $\varphi'$ and $\psi$ (Figure \ref{figura11}) ;
let $[\frac{1}{4},z_1]$ and $ [d_0,d_1]$ be the intervals such that $ \varphi '(z) \in [\frac{1}{2},1] \Leftrightarrow z \in [\frac{1}{4},z_1] $ and $ \psi(x^2+y^2) \in [2,4] \Leftrightarrow x^2+y^2 \in [d_0,d_1] $.\\
  For each $z \in [\frac{1}{4},z_1]$ there exists a real value  $ d=d(z) \in [d_0,d_1]$ with $\psi(d(z))=\frac{2}{\varphi'(z)}$; it holds that $\{x^2+y^2=d(z)\} \times {\{z\}} \subset S_f$. Moreover, there exists a unique value $d=\frac{1}{16}$ for the cases $z=\frac{1}{4}$ or $z=z_1$ and there exist \emph{two} distinct values $d(z)<D(z)$ for each $z \in (\frac{1}{4},z_1)$ that verify this property (they are symmetrical with respect to $\frac{1}{16}$). This situation can be regarded as that while $z$ 'runs through' the interval $[\frac{1}{4},z_1]$, the critical set gets foliated as the union of two circles. These circles are defined by the intersection of the boundaries of two balls centered at $(0,0,z)$, with radii $d(z)$ and $D(z)$, with 'horizontal' hyperplanes defined by constant values of $z$ as Figure \ref{nigual3} shows (observe that these circles are spheres of codimension 2 in $\T^3$).\\
   We have shown then that $S_f=\bigcup_{z \in [\frac{1}{4},z_1]}(\{x^2+y^2=d(z)\} \cup \{x^2+y^2=D(z)\}) \times \{z\} $ or, in the same fashion, $S_f=\bigcup_{z \in [\frac{1}{4},z_1]}(\partial B_{(\tilde{0},\sqrt{d(z)})} \cup \partial B_{(\tilde{0},\sqrt{D(z)})} ) \times \{z\} $. Finally, since $z=\frac{1}{4}$ and $z=z_1$ have the same values of $d(z)=D(z)=\frac{1}{16}$ it is possible to identify the boundary manifolds $\{x^2+y^2=\frac{1}{16}\}\times\{\frac{1}{4}\}$ and $\{x^2+y^2=\frac{1}{16}\}\times\{z_1\}$ to obtain \textbf{$S_f \simeq S^{1} \times S^1$}. Observe that the extremal values for the radii of the balls are obtained when $z=\frac{1}{4}+\frac{\delta}{8}$ (they appear in Figure \ref{nigual3} as $B_0(x_\delta)$ and $B_1(x_\delta)$ respectively). \\

\begin{figure}[ht]
\begin{center}
\includegraphics[scale=0.25]{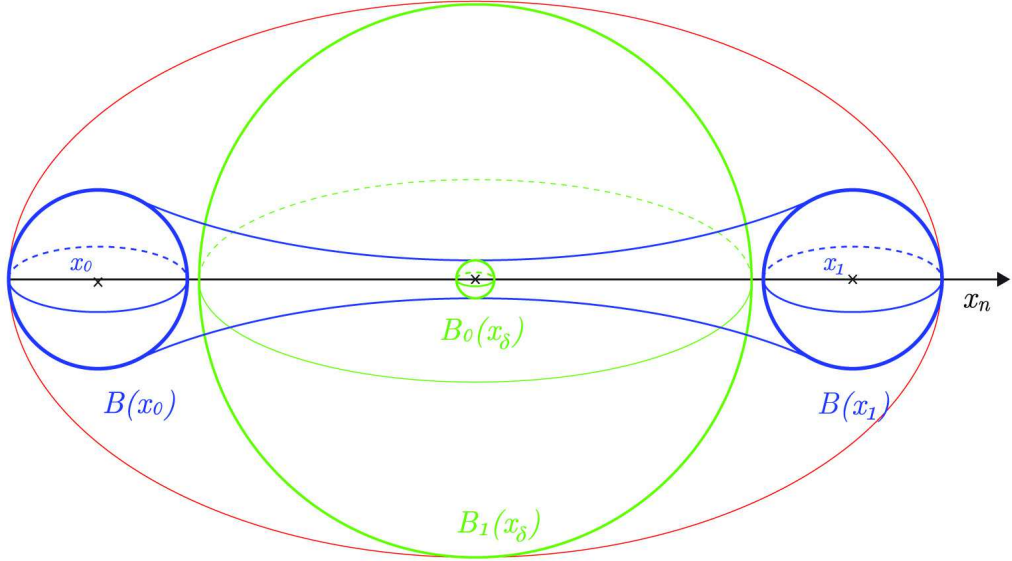}
\caption{Frontal view of $S_f$ for $n \geq 3$, $x_0=\frac{1}{4}$, $x_1=c$. }\label{Sf}
\end{center}
\end{figure}

Analyzing the \underline{general case} as the previous ones we deduce:\\
 $S_f=\{x \in \T^{n} / 2-\varphi'.\psi =0\rbrace$ and $\psi \leq 4$ so it has to be $\varphi' \geq \frac{1}{2}$. In consequence, for all critical points $x=(x_1,...,x_n) \in S_f$ it hold $ \varphi'(x_n) \geq \frac{1}{2}$ and $\psi \in [2,4]$. Turn once more to the graphs of $\varphi'$ and $\psi$ (Figure \ref{figura11}).\\
  Let $[\frac{1}{4},c]$ be the interval where $ \varphi '(x_n) \in [\frac{1}{2},1] \Leftrightarrow x_n \in [\frac{1}{4},c] $.\\
Let $ [d_0,d_1]$ be the interval where $ \psi(d) \in [2,4] \Leftrightarrow d \in [d_0,d_1] $.\\
For each  $x_n \in [\frac{1}{4},c]$ there exists a real value $ d=d(x_n) \in [d_0,d_1]$ such that $\partial B_{(\tilde{0},\sqrt{d(x_n)})} \times {\{x_n\}} \subset S_f$. And just like before, for every $  x_n \in (\frac{1}{4},c) $ there exist \emph{two} distinct values $d(x_n)<D(x_n)$ that satisfy this condition except for a unique value $d(x_n)=D(x_n)=\frac{1}{16}$ at both cases $x_n=\frac{1}{4}$ and $x_n=c$. What happens is that while $x_n$ 'runs through' the interval $[\frac{1}{4},c]$, the critical set is foliated by the boundary of two balls of codimension 1, which are spheres of codimension 2 in $\T^n$. Once again it allows writing the critical set as $S_f=\bigcup_{x_n \in [\frac{1}{4},c]}(\partial B_{(\tilde{0},\sqrt{d(x_n)})} \cup \partial B_{(\tilde{0},\sqrt{D(x_n)})}) \times \{x_n\} $. Moreover, since the balls at $x_n=\frac{1}{4}$ and $x_n=c$ are equal (appearing in Figure \ref{Sf} as $B(x_0)$ and $B(x_1)$ respectively) we can identify them to obtain \textbf{$$S_f \simeq S^{n-2} \times S^1$$} which is our claim. Again, the extremal values for the radii of these balls are obtained at the \textit{equator},  when $x_n=\frac{1}{4}+\frac{\delta}{8}$ (they appear in Figure \ref{Sf} as $B_0(x_\delta)$ and $B_1(x_\delta)$ respectively). $\Box$

\subsubsection{Classification of the critical points}\label{folds}

We pay attention now to a special type of singularity called \textit{fold points} which are the most frequent singularities to appear on any map, corresponding to symbol $(1,0)$ in the Thom-Boardman classification. In particular, we pay special attention to these since they are the only type of singularity that will appear explicitly along this work, and most important of all because near fold points open sets collapse (as Lemma \ref{abiertofold} shows), and this can of course prevent robust transitivity as we will show happens for the map $f$ studied in the previous sections.\\ In the classical development of singularity theory, fold points (as well as all of the others) are defined by their normal forms and then its properties are established. We will not do so, but follow the path presented in \cite{gg} and give the following definition (we won't go on about their normal form until we need it later on in the article):

\begin{defi}\label{defifold}
  We say that $x \in S_f$ is a \textbf{fold singularity}  if $ ker(D_xf) \transv T_xS_f$.
\end{defi}

We denote by $F_f$ the set of fold points of $f$; i.e. $F_f=\{ x \in S_f / x \mbox{ is a fold point} \}$.

  \begin{thm}\label{clasificaciondelospuntoscriticos}
    Let $f$ be the map defined by Equation (\ref{mapaefe}). Then, every critical point in $S_f$ except for those that satisfy $x_n=\frac{1}{4}+\frac{\delta}{8}$ is a fold point. Or equivalently, $S_f \setminus F_f \subset \{ x_n= \frac{1}{4}+\frac{\delta}{8} \}$.
  \end{thm}
  \pf
  Recall observation 2 in Remark \ref{rkcritset}. It provides $ker(D_x f)=\langle\{\vec{e_n}\}\rangle$ for all $x \in S_f $. Since $dim(T_x S_f)=n-1$, it suffices that $\{(0,0,...,0,1)\} \nsubseteq T_x S_f$ for $x$ to be a fold point. According to Lemma \ref{lema 2} we know that the tangent space is given by
  $T_x S_f = T_x h^{-1}(\mathcal{R}_{n-1})=ker(D_x h)$. On the proof of the Lemma it was shown that if $\varphi^{''} \neq 0 $ then $(0,0,...,1) \in Im(D_x h)$ hence $ (0,0,...,1) \notin ker(D_x h)$ which implies $x \in F_f$. Therefore, whenever a critical point $x$ is not a fold point then $\varphi^{''} = 0$ which means $ x_n=\frac{1}{4}+\frac{\delta}{8}$. $\Box$
 \begin{clly}
   $F_f$ is open and has total Lebesgue measure in $S_f$.$\Box$
 \end{clly}
\begin{rk}
  What Theorem \ref{clasificaciondelospuntoscriticos} shows is that all points in $S_f$ except for those at the 'equators' are fold points. Geometrically it is easy to see why, since being $ker(D_xf)$ generated by $\vec{e_n}$ for all $x \in S_f$, the equators are the only lines where the transversality condition fails to hold.
\end{rk}

\section{A singular endomorphism $H$ of $\T^n$}

Now that we have thoroughly studied the critical set and dynamical behavior of the map $f$, we turn to the construction of the objective map of this work. That is, a map $H$ that is $C^2$ robustly transitive but not $C^1$ robustly transitive which we obtain through a perturbation of $f$. The idea is to see that, given $\ep >0$, around $p$ there exists an open set of $\T^n$ whose image under some $g$ in $\mathcal{U}_{(H,\varepsilon)} \in C^1$ is meager and invariant which prevents robust transitivity. Finally we will prove that $H$ is $C^2$ robustly transitive making use of the unstable cones it inherits from $f$.

\subsection{Construction of $H$.}\label{H}

Recall $p=(\frac{1}{4},0,..,0,\frac{1}{4}) \in S_{f} $. We start showing that around $p$ all critical points have their last coordinate as an implicit function of the others. We show as well that every first and second order derivatives of this implicit function at $p$ are equal to zero.

\begin{thm}\label{implicita}
 There exist $U \subset \R^{n-1}$ open neighborhood of $\tilde{p}=(\frac{1}{4},0,...,0)$; $ V \subset \R$ open neighborhood of $\frac{1}{4}$ and a smooth function $ \phi:U \rightarrow V$ with $\phi(\tilde{p})=\frac{1}{4}$ such that for all $ x=(x_1,...,x_n) \in S_f \cap U\times V $ it is $ x_n=\phi(x_1,...,x_{n-1})$.
\end{thm}
\pf
Recall $S_{f}=\lbrace x \in \T^n /2-\varphi'.\psi =0\rbrace $. Derivate over $S_f$ along $x_n$ to obtain $ \frac{\partial}{\partial x_{n}}( 2-\varphi'.\psi )=-\varphi'' .\psi$. The evaluation at $p$ is not zero since $\varphi''(\frac{1}{4})\neq 0$.
Apply the Implicit Function Theorem to obtain an open neighborhood $U$ of $\tilde{p}=(\frac{1}{4},0,...,0)$, $ V$ an open neighborhood of $ \frac{1}{4}$ and a smooth function $ \phi:U \rightarrow V $ satisfying $ \phi(\tilde{p})=\frac{1}{4}$ and such that for all  $ x \in S_f \cap U \times V, x=(\tilde{x},\phi(\tilde{x}))$.$\Box$ \\

We prove now that all first and second order derivatives of $\phi$ are null at $\tilde{p}$.

\begin{lema}\label{valecero}  $\forall k, k' \leq n-1: \frac{\partial}{\partial x_{k}}\phi(\tilde{p})=\frac{{\partial}^2}{\partial x_{k}.x_{k'}}\phi(\tilde{p})=0$.
\end{lema}
\pf
For all $  x \in S_f \cap U \times V$, it holds $ 2-\varphi '(\phi).\psi=0 $. Let $ k\leq n-1$ and derivate along $x_k$; it holds $ \frac{\partial}{\partial x_k}[2-\varphi'(\phi).\psi]=0 $, then $$ \forall k\leq n-1: -\varphi^{''}(\phi).\frac{\partial}{\partial x_{k}}\phi.\psi-2.\varphi'(\phi).x_{k}.\psi '=0 \therefore \frac{\partial}{\partial x_{k}}\phi=-\frac{2x_{k}.\varphi'(\phi).\psi'}{\varphi^{''}(\phi).\psi}. $$ \\ Let $ k'\leq n-1$ and derivate now along $x_{k'}$: $$ \frac{{\partial}^2}{\partial x_{k}.x_{k'}}\phi= \left[ \frac{\partial}{\partial x_{k'}}\left(\frac{-2x_k.\varphi'(\phi)}{\varphi^{''}(\phi).\psi}\right)\right].\psi'+  \left(\frac{-2x_k.\varphi'(\phi)}{\varphi^{''}(\phi).\psi}\right) .\left( 2x_{k'}  \psi^{''}\right).$$\\ It is clear that both of them are well defined in $U$ and since $\psi^{'}( \frac{1}{16})=\psi^{''}( \frac{1}{16})=0$ they become zero at $\tilde{p}$.$\Box$

\begin{clly}
  $\forall\varepsilon' >0, \exists U' \subset U$ open neighborhood of $\tilde{p}$ in $\R^{n-1} $ such that \\$\forall x \in U'; \forall k, k' \leq n-1; $ $$ max\left\{\left|\phi(x)-\frac{1}{4}\right|,\left|\frac{\partial}{\partial x_{k}}\phi(x)\right|,\left|\frac{{\partial}^2}{\partial x_{k}.x_{k'}}\phi(x)\right|\right\}<\varepsilon'. $$
\end{clly}
\pf Apply Lemma \ref{valecero} and definition of continuity.$\Box$
\begin{clly}
  $\forall\varepsilon' >0, \exists U' \subset U$ open neighborhood of $\tilde{p}$ in $\R^{n-1}$ such that \\$ \forall x \in U'; \forall k \leq n-1;$ $$ \left| \phi(x) -\frac{1}{4} \right| < \varepsilon'\left\|x-\tilde{p}\right\|^2 \mbox{ and } \left| \frac{\partial}{\partial x_k} \phi(x) \right| < \varepsilon'\left\|x-\tilde{p}\right\|. $$
\end{clly}
\pf Apply Taylor's Theorem to $\phi$ at $\tilde{p}$ together with Lemma \ref{valecero}.$\Box$
\\

We turn our attention now to the image $f(S_f)$ of the critical set under $f$. Here again, near $f(p)$ all points in $f(S_f)$ have their last coordinate as a function $\Phi$ of the previous ones, and all first and second order derivatives of this function $\Phi$ are null at $\widetilde{f(p)}$. Observe that since $A(p)=f(p)=(\tilde{0},\frac{1}{2})$ then $ \widetilde{A(p)}=\widetilde{f(p)}=\tilde{0}$. \\
To make the reading easier we denote through the rest of the section $\phi(\tilde{x})$ as $\tilde{\phi}$ for all $\tilde{x} \in U$. Theorem \ref{implicita} grants that for all points $x$ in $S_f \cap U \times V$ it holds that $$f_{|S_f \cap U\times V}(x_1,...,x_n)=(8x_1,2x_2,...,2x_{n-1},2\tilde{\phi}-\varphi(\tilde{\phi}).\psi).$$
Observe now that since the matrix $\widehat{A}$ at the start of Section 3 is invertible allows the definition of a map $ A^{-1}:\R^n \rightarrow \R^n /A^{-1}(x_1,...,x_n)= (\frac{x_1}{8},\frac{x_2}{2},...,\frac{x_n}{2}) $. This map, fairly being not defined on $\T^n$, is useful for finding the analytic expression of the function $\Phi$ we seek since a change of variables $z:=\widetilde{A(x)}$ enables writing $$f_{|{S_f \cap U\times V}}(x)=(z,2\tilde{\phi}(A^{-1}((z,0)))-\varphi(\tilde{\phi}(A^{-1}((z,0)))).\psi(A^{-1}((z,0)))).$$ Define  $$\Phi:\R^{n-1}\rightarrow \R / \Phi(z)=[2.(\tilde{\phi}\circ A^{-1})-(\varphi \circ \tilde{\phi}\circ A^{-1}).(\psi \circ A^{-1})](z,0)$$ to have an explicit analytic expression for the last coordinate as a function of the others for all points near $f(p)$ in the set $f(S_f \cap U \times V)$.

\begin{rk}\label{ultimacoordenada}
  For all  $ x \in f(S_f \cap U \times V)$ the equality $x_n=\Phi(\tilde{x})$ holds.
\end{rk}

\begin{lema}
  $\forall k, k' \leq n-1: \frac{\partial}{\partial x_k}\Phi(\tilde{0})=\frac{\partial^2}{\partial x_k.x_{k'}}\Phi(\tilde{0})=0$.
\end{lema}
\pf It suffices to notice that $A_{|\{x_n=0\}}^{-1}$ is linear. Apply the chain rule when derivating $\Phi$ and Lemma \ref{valecero} to obtain the proof to the claim. $\Box$\\

Moreover, considering that $A_{|\{x_n=0\}}^{-1}$ contracts and that $|\varphi'|\leq1$, the following corollaries are in order:

\begin{clly}\label{A1}
  $\forall\varepsilon' >0, \exists W $ open neighborhood of $\tilde{0}$ in $\R^{n-1}$ such that \\ $\forall x \in W;  \forall k, k' \leq n-1; $ $$ max\left\{\left|\Phi(x)-\frac{1}{2}\right|,\left|\frac{\partial}{\partial x_{k}}\Phi(x)\right|,\left|\frac{{\partial}^2}{\partial x_{k}.x_{k'}}\Phi(x)\right|\right\}<\varepsilon'.\Box$$
\end{clly}

\begin{clly}\label{A2}
  $\forall\varepsilon' >0, \exists W $ open neighborhood of $\tilde{0}$ in $\R^{n-1}$ such that \\ $\forall x \in W;  \forall k \leq n-1;$ $$ \left| \Phi(x) -\frac{1}{2} \right| < \varepsilon'\|x\|^2 \mbox{ and } \left| \frac{\partial}{\partial x_k} \Phi(x) \right| < \varepsilon'\|x\|.\Box$$
\end{clly}
\begin{figure}[ht]
\begin{center}
\includegraphics[scale=0.35]{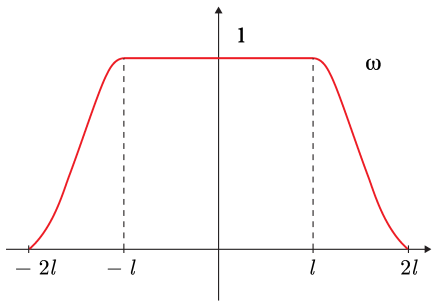}
\caption{Graph of $\omega$}\label{figura22}
\end{center}
\end{figure}
From now on we proceed to build a diffeomorphism $F$, $C^2$ close to the identity map, that we will compose with $f$ to obtain $H$. Essentially, what $F$ will do is a 'flattening' of $f(S_f)$ near $f(p)$ and be the identity map away from it. \\ To begin with, given $\ep >0$, choose $W$ where corollaries \ref{A1} and \ref{A2} hold for an $\varepsilon' <\frac{\varepsilon}{50}$. Next, fix a parameter $b \in \R $ such that $0 < b <\frac{1}{4}$ satisfying  that the ball centered at $\tilde{0}$ and radio $2b$ has closure contained in $W$. Then, given $l \in \R$ such that $0<l < b$, define an auxiliary smooth function $\omega: \R \rightarrow\R$ (as Figure \ref{figura22} shows) satisfying:
\begin{itemize}
  \item $\omega_{| B_{(0,l)}}=1$, $\omega_{| B_{(0,2l)}^{c}}=0$ and $\|\omega \|=1$,
  \item $\|\omega'\|<\frac{2}{l}$ and $\|\omega''\|<\frac{8}{l^2}$.
\end{itemize}

\vspace {0.2cm}
Finally, define $u : \R^{n-1}\rightarrow \R$ $/u(\tilde{x})=\omega(\|\tilde{x}\|).(\Phi(\tilde{x})-\frac{1}{2})+\frac{1}{2}$; \\observe that $u_{| B_{(\tilde{0},l)}}=\Phi$ and $u_{| B_{(\tilde{0},2l)}^{c}}=\frac{1}{2}$.

\begin{lema}\label{phichato}
  $$ \forall \tilde{x} \in \R^{n-1}, \forall k, k' < n: max\left\{ \left|u(\tilde{x})-\frac{1}{2}\right|,\left|\frac{\partial}{\partial x_k}u(\tilde{x})\right|,\left|\frac{\partial^2}{\partial x_k \partial x_{k'}}u(\tilde{x})\right| \right\} <\varepsilon. $$
\end{lema}

\pf
It suffices to calculate the extreme values in $\overline{B_{(\tilde{0},2l)}\setminus B_{(\tilde{0},l)}}$ since in the bounded component of the complement is less or equal than $\varepsilon'$ (controlled by $\Phi$) and in the unbounded component is null. Moving on to calculations for $\tilde{x}$ inside the annulus, we write  $\omega(\|\tilde{x}\|)$ as $\omega$ in the sake of simplicity. Compute derivatives and apply Corollary \ref{A2}:\\
1) $$ \left|u(\tilde{x})-\frac{1}{2}\right| = \left|\omega.\left(\Phi(\tilde{x})-\frac{1}{2}\right) \right| \leq \left|\Phi(\tilde{x})-\frac{1}{2}\right| \leq \varepsilon' <\varepsilon . $$\\
2)$$\left|\frac{\partial}{\partial x_k}u(\tilde{x})\right|=\left| \omega'.\frac{(\tilde{x})_k}{\|\tilde{x}\|}.\left(\Phi(\tilde{x})-\frac{1}{2}\right)+\omega.\frac{\partial}{\partial x_k} \Phi(\tilde{x}) \right| \leq \frac{2}{l}4l^2\varepsilon'  +2l\varepsilon'  < 10l\varepsilon'<\varepsilon. $$\\
3) $$\left|\frac{\partial^2}{\partial x_k \partial x_{k'}}u(\tilde{x})\right|= |\omega^{''}.\frac{(\tilde{x})_{k'}}{\|\tilde{x}\|}.\frac{(\tilde{x})_k}{\|\tilde{x}\|}.\left(\Phi(\tilde{x})-\frac{1}{2}\right)+\omega^{'}.\frac{-(\tilde{x})_{k}.(\tilde{x})_{k'}}{\|\tilde{x}\|^3}.\left(\Phi(\tilde{x})-\frac{1}{2}\right)+$$ $$+ \omega'.\frac{(\tilde{x})_k}{\|\tilde{x}\|}.\frac{\partial}{\partial x_{k'}} \Phi(\tilde{x})+ \omega'.\frac{(\tilde{x}-\tilde{f}(p))_{k'}}{\|\tilde{x}-\tilde{f}(p)\|}.\frac{\partial}{\partial x_k}\Phi(\tilde{x})+\omega.\frac{\partial^2}{\partial x_k \partial x_{k'} } \Phi(\tilde{x})| \leq $$ $$ \leq  \frac{8}{l^2}4l^2\varepsilon'+\frac{2}{l}\frac{1}{l}4l^2\varepsilon'+2(\frac{2}{l}2l\varepsilon')+\varepsilon'<50\varepsilon'<\varepsilon. \Box $$
\begin{clly}\label{coroFachata}
   Given $\ep >0$, there exists a diffeomorphism $F$ of $\T^n$ such that $d_{C^2}(F,Id) < \varepsilon$ and exists $l>0 $ such that  $\forall x \in f(S_f) \cap B_{(f(p),l)}:  F(x)=(\tilde{x},\frac{1}{2}). $
\end{clly}
\pf
Define explicitly $F: \T^n \rightarrow \T^n / F(x)=(\tilde{x},x_n-u(\tilde{x})+\frac{1}{2})$. \\
Then,
\begin{itemize}
\item $F$ is a diffeomorphism of  $\T^n$ since $ det(D_xF)=1, \forall x \in \T^n$.
  \item $ \forall x \in \T^n, \|(F-Id)(x) \|=|u(\tilde{x})-\frac{1}{2}| $. By Lemma  \ref{phichato},  $d_{C^2}(F,Id)<\varepsilon$ .
    \item $F_{|B_{(\tilde{0},2l)}^c \times \R}(x)=Id$ therefore $F_{|B_{(f(p),2l)}^c}(x)=Id$.
  \item $F_{|B_{(\tilde{0},l)} \times \R}(x)=(\tilde{x},x_n + \frac{1}{2}-\Phi(\tilde{x}))$.
  \item By remark \ref{ultimacoordenada}, $ x_n=\Phi(\tilde{x})$ for all points $x \in f(S_f \cap U \times V)$. Then, the restriction $F_{|B_{(f(p),l)}\cap f(S_f \cap U \times V) } \subset\{ x_n= \frac{1}{2}\}$.$\Box$
\end{itemize}
\vspace{0,4cm}
We are now in condition to define a map $H$ that has the properties claimed in the title of the article:
 $$H:\T^n \rightarrow \T^n / H(x)=(F \circ f)(x).$$

\begin{rk}\label{remarkH}
  It is straightforward seeing that the following hold:
\begin{itemize}
\item $d_{C^2} (H,f)< \varepsilon$.
  \item  $H(p)=f(p)$.
  \item $S_H=S_f$ hence $p \in S_H$.
  \item $H(B_{(p,r)})\cap B_{(p,r)} = \emptyset$ because of how $f$ and $r$ were defined at the beginning.
  \item $H(x)=A(x)$, $ \forall x \in \T^n \setminus (B_{(p,r)} \cup A^{-1}(B_{(A(p),r)}))$ because of how $H$ is defined.
  \item There exists $ a' \leq a $ such that $\overline{D_{x}H(C^u_{a^{'}}(x))} \setminus \{ (0,0) \} \subset C^u_{a^{'}}(H(x)) , \forall x \in \T^n$ by Corollary \ref{coroconos}. Observe that the parameter $\varepsilon'$ of Corollaries \ref{A1} and \ref{A2} as well as the parameter $l$ from Lemma \ref{phichato} can all be shrunk if it was desired to push $H$ closer to $f$.\item $H(S_H)= F[f(S_f)]  \subset \{x_n=\frac{1}{2}\}$, $\forall x \in B_{(f(p),l)}$ by Corollary  \ref{coroFachata}.
\end{itemize}
\end{rk}

We proceed in the following sections to prove that $H$ is an endomorphism of $\T^n$ that is $C^2$ robustly transitive but not $C^1$ robustly transitive.

\subsection{$H$ is not $C^1$ robustly transitive}

We begin the section showing that every map with a fold singularity is $C^1$ close to a map whose critical set has nonempty interior.\\

 Remember that $x \in S_f$ is a fold singularity if $ ker(D_xf) \transv T_xS_f$ and that the set $F_f=\{ x \in S_f / x \mbox{ is a fold point} \}$. In \textit{\cite{gg},chapter 3, Theorem 4.5} it is shown that for every $x \in F_f$, there exists a change of coordinates in $\R^n$ such that $f$ takes its \textbf{normal form }. Explicitly, if $f:M \rightarrow M / f(x)=y$ and $x \in F_f$ then there exist $ U, V$ open sets containing $x$ and $y$, there exists $N$ open neighborhood of $0 \in \R^m$ and diffeomorphisms $g:N \rightarrow U$ and $h:V \rightarrow N$ such that for all $ (x_1,...,x_n) \in N$, $( h \circ f_{|U} \circ g)(x_1, ..., x_n)=(x_1,...,x_{n-1},x_n^2)$.

\begin{figure}[ht]
\begin{center}
\includegraphics[scale=0.35]{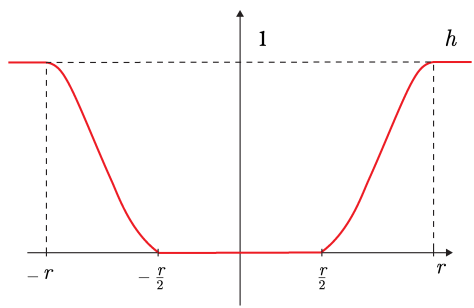}
\caption{Graph of $\mu$}\label{uno}
\end{center}
\end{figure}

\begin{lema}\label{abiertofold}
 If $f \in C^1(M,M)$ and $x \in F_f$ then for every $ \varepsilon >0$ and for every $W$ open neighborhood of $x$ in $M$ there exist a map $g$ in a neighborhood $ \mathcal{U}_{(f,\varepsilon)} \in C^1$ and $V$ an open neighborhood of $x$ in $M$ such that $V \subset S_g, g(V)=f(S_f \cap V)$ and $g(x)=f(x),\forall x \notin W$.
\end{lema}
\pf
Without loss of generality we take $x=0$ and $f$ in \emph{normal form}. \\ Let $\varepsilon >0$ and $W$ an open neighborhood of $x$ be given and choose $r>0$ such that $B_{(0,r)} \subset W$ and that $ 7r<\varepsilon$ and $r^2<\varepsilon$.\\ Take a bump function $\mu:\R \rightarrow \R$ class $C^1$ as shown in Figure \ref{uno} with $|\mu'|<\frac{4}{r}$ and define a map $g:\R^n\rightarrow \R^n/ g(x_1,..,x_n)=(x_1,...,x_{n-1},x_n^2.\mu(x_n))$. A straightforward calculation shows that $(f-g)(x_1,..,x_n)=(0,...,0,x_n^2.(1-\mu(x_n)))$ in $B_{(0,r)}$ and zero in $B_{(0,r)}^c$. Another straightforward calculation provides:

\begin{itemize}
  \item $|f-g|\leq |x_n^2|.|1-\mu(x_n)|<r^2<\varepsilon$.
  \item $|\frac{\partial}{\partial x_n}(f-g)|\leq |2x_n(1-\mu(x_n))-x_n^2.\mu'(x_n)|<2r+r^2.\frac{4}{r}<7r<\varepsilon$.
\end{itemize}
Define $V=B_{(0,\frac{r}{2})}$ so then $ g_{|V}=(x_1,..,x_{n-1},0)=f_{|S_f}=f_{|S_f \bigcap V}$. The remaining assertion follows from the fact that $W^c \subset B_{(0,r)}^c$ .$\Box$

\begin{rk}
  Lemma \ref{abiertofold} does not hold in a $\mathcal{U}_{(f,\varepsilon)} \in C^2$ neighborhood of $f$. It is sufficient to notice that the Hessian matrix of $f$ has an element $2$ in the last entry. We provide with a more general proof for this remark in Theorem \ref{nopasac2} where we don't assume the singularity to be in normal form nor perform computations of derivatives.
\end{rk}

Now that the key feature about fold points has been stated, we move on closing the subsection with the proof of the claim on its title.

\begin{lema}
  $p$ is a fold point for $H$.
\end{lema}
\pf
Let $h$ be the map defined by Equation \ref{mapah},  then Remark \ref{remarkH} provides
$T_p S_H = T_p S_f = T_p h^{-1} (\mathcal{R}_{n-1}) = ker (D_p h)$ and $ ker (D_p H)= ker(D_p f)$.
    Apply Lemma \ref{lema 2} and remark 2 in Equation \ref{ecuacion1} to obtain $T_p S_H \oplus ker (D_p H) = \R^n $ which is the condition for $ p \in F_H$ . $\Box$

\begin{figure}[ht]
\begin{center}
\includegraphics[scale=0.16]{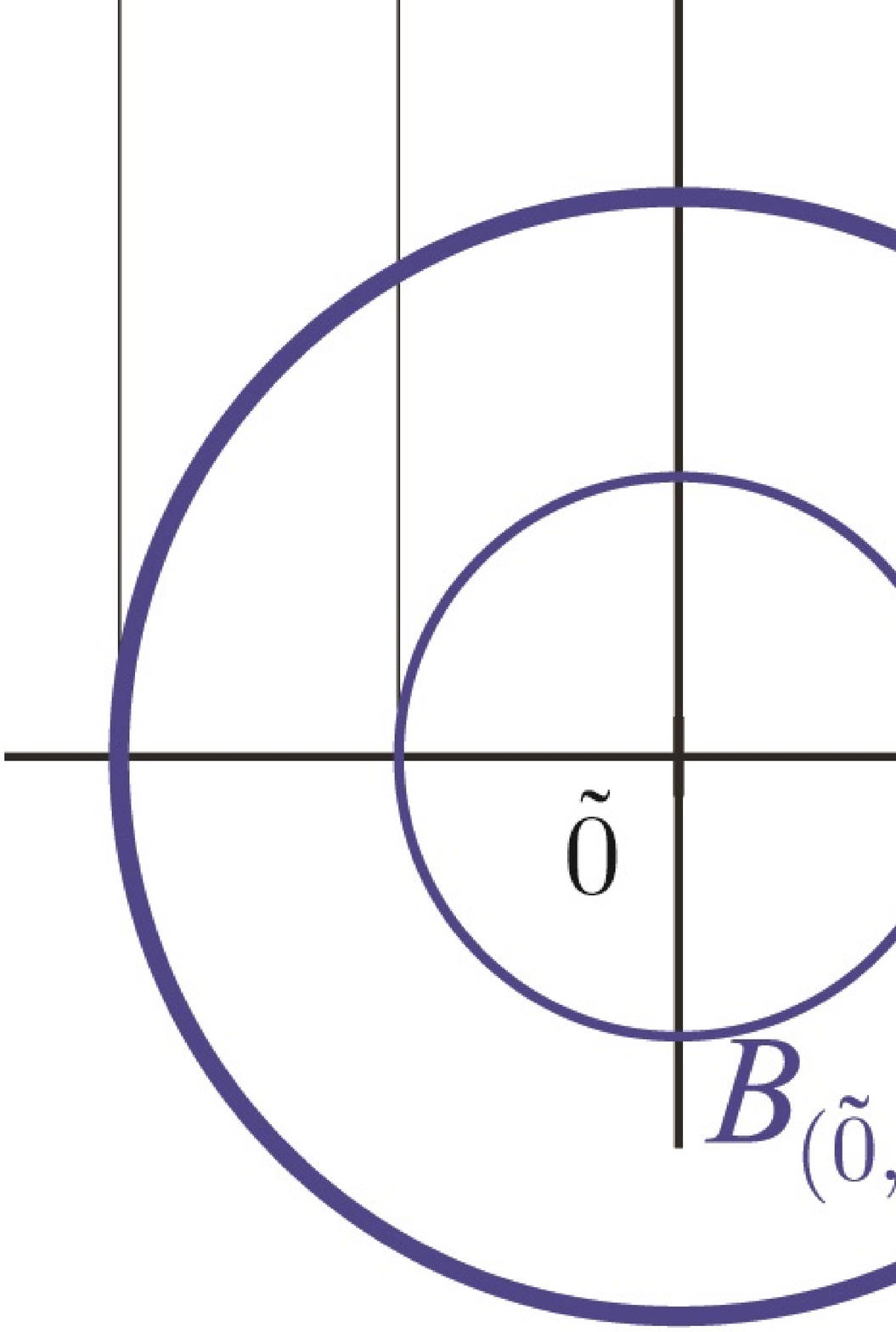}
\caption{$H$ 'flattens' the image of $S_f$ around $f(p)$}\label{Haplasta}
\end{center}
\end{figure}

\begin{thm}
  $H$ is not $C^1$ robustly transitive.
\end{thm} \pf (see Figure \ref{Haplasta})\\
Fix any $\varepsilon >0$. Since $p \in F_H$, Lemma \ref{abiertofold} provides a function $g $ in $\mathcal{U}_{(H, \ep)} \in C^1$ and an open set $W$ in $\T^n$ containing $p$ such that: $W \subset S_g$, $g(W)=H(S_H \cap W)$ for all $ x \in W$ and $g=H$ in $B_{(p,l)}^c$.
Abusing language in Corollary \ref{coroFachata}, call $W \cap B_{(p, l)}$ as $W$, then $g(W)=H(S_H \cap W) \subset H(S_H) \subset \{x_n=\frac{1}{2}\}$. Since $H=A$ in $B_{(p,r )}^c$, $f(p)\neq p$ and $|g(p)-f(p)|<\varepsilon$ then, previously reducing $W$ again if needed, it holds that  $g^2(W)=g(g(W))=H(g(W))=A(g(W)) \subset \{x_n=0\}$ . On the other hand, $\{x_n=0\}$ is invariant under $g$ (since 'far' from $p$, $g$ is equal to $A$) then $ g^m(W) \subset \{ x \in \T^n / x_n=0  \} ,\forall m \geq 2$ which implies $g$ is not transitive. Being $g$ not transitive, $H$ is not $C^1$ robustly transitive. $\Box$

\subsection{$H$ is $C^2$ robustly transitive}

We start the final subsection with a theorem showing that in the $C^2$ topology it is possible to find a neighborhood of $f$ where no map in there collapses open sets; opposed to what happens in the $C^1$ topology stated in Lemma \ref{abiertofold}:

\begin{thm}\label{nopasac2}
   Let $f$ be the map defined by Equation \ref{mapaefe}. There exists a $C^2$ neighborhood $\mathcal{U}_f$ of $f$ such that for all $ g \in \mathcal{U}_f$ and all $V$ open in $\T^n$, the interior $int(g(V))$ is not empty.
\end{thm}
\pf
Observe that Lemmas \ref{lema 1} and \ref{lema 2} hold simultaneously for some open set $\mathcal{U}_f \in C^2$ containing $f$, therefore for all $ g \in \mathcal{U}_f$, $S_g$ is a submanifold of $\T^n$ of codimension 1. Given an arbitrary $g$ in $\mathcal{U}_f$ and any $ V$ open in $\T^n$ there exists $U \subset V$ open such that $U \cap S_g = \emptyset$. Since $g$ is a local diffeomorphism in the interior of $V \setminus S_g$ it carries $U$ to an open set. Hence, the interior of $g(V)$ is nonempty. $\Box$\\

We point out now the key properties that $f$ satisfies which are needed to finish the proof:
\begin{itemize}
 \item  1) $f$ admits unstable cones and this is a $C^1$-stable property according to Corollary \ref{coroconos} (is worth to mention that the parameter $a$ can be chosen as small as desired in order to obtain a bounded prefixed difference between diameter and length, let's say 'less than $\ep$'. Hence diameter and length can be identified without risk).
  \item 2) $f_{| \T^n \setminus B_{(p,r)}}$ expands due to Observation \ref{obs1}.
\end{itemize}

Define then a $C^1$ neighborhood of $f$ such that all of its elements admit unstable cones (item 1 holds). Then, reduce it until  for all maps $ \hat{f}$ in there hold that $\hat{f}_{| \T^n \setminus B_{(p,\frac{3}{2}r)}}$ expands (item 2 holds) and reduce it again to a neighborhood $\mathcal{{U}}_f \in C^2$ of $f$ where Theorem \ref{nopasac2} also holds for all $ \hat{f} \in \mathcal{U}_f$. Also, since $H=F \circ f$ and $d_{2}(F,Id)<\varepsilon$ we can take $H \in \mathcal{U}_f$ and $\mathcal{U}_H \in C^2$ such that $\mathcal{U}_H \subset \mathcal{U}_f$. \\

In what comes next we choose an arbitrary $g \in \mathcal{U}_H$. If we prove $g$ is transitive then we will have proved that $H$ is $C^2$ robustly transitive. The strategy relies in noticing that for all open sets $V$ there exist within $V$ curves travelling inside $g$'s unstable cones which escape the perturbation region without coming back inside in the future. Being in a $C^2$ neighborhood where Theorem  \ref{nopasac2} holds, then an open subset of $V$ also escapes the perturbation region. Finally, being $g$ \textit{$\ep$-close} to $A$ far away from $p$, then $g$ will be transitive. \\ We prove a series of lemmas that lead to the formal proof of these claims.

\begin{lema}\label{sucesiondecurvas}
  Let $V \subset \T^n$ open. There exist a point $ y \in V$ and $ n_0 \in \N$ such that $g^k(y) \in int(g^k(V))$ and $g^k(y) \in \T^n \setminus B_{(p,2r)}, \forall k \geq n_0$.
\end{lema} \pf
Let $\beta \subset int(V)$ be a curve parametrized over an interval $I$ such that $\forall t \in I, \beta'(t) \in C^u_a(\beta(t))$. Condition \ref{3} on the differential of $f$ and the invariance of the unstable cones grant $diam(g^n(\beta))> 6^n.diam(\beta)$. Therefore it exists $n_0$ such that $diam(g^{n_0}(\beta)) > 9r$. Take a curve $\alpha \subset g^{n_0}(\beta) \setminus B_{(p,2r)}$ with $diam(\alpha) \geq 2r$.
Build now a succession of curves $\{\alpha_k\}$ in the following fashion: take $\alpha_1=\alpha$, since  $\alpha'(t) \in C^u_a(\alpha(t))$ then $diam(g(\alpha_1))> 12r$ .
Being $diam(B_{(p,2r)})=4r$,  there exists $ \alpha_2 \subset g(\alpha_1)$ such that $ \alpha_2 \subset \T^n \setminus B_{(p,2r)}$ and $diam (\alpha_2) \geq 2r$. Proceed inductively according to this algorithm to find a family of curves $\{ \alpha_k \}$ satisfying for all $ k \in \N $:
\begin{itemize}
 \item $g(\alpha_{k}) \supset \alpha_{k+1}$,
\item $\alpha_k \subset \T^n \setminus B_{(p,2r)} $,
\item $\alpha'_k(t) \in C^u_a(\alpha_k(t))$.

\end{itemize}

  Define afterwards a family of curves $\gamma_k$ such that $\forall k: \gamma_k \subset \alpha_1$ and $g^k(\gamma_k)=\alpha_k$.
We have  $g^{k+1}(\gamma_k)=g(g^k(\gamma_k))=g(\alpha_k) \supset \alpha_{k+1}$, then $g^{k+1}(\gamma_k) \supset g^{k+1}(\gamma_{k+1})$ which imply $\gamma_k \supset \gamma_{k+1}$. Apply Cantor's Intersection Theorem to obtain a point $ x$ such that $\bigcap_{k \in \N} \overline{\gamma_k} =\{x\}$. It holds for all $k \in \N$ that $g^k(x) \in \T^n \setminus B_{(p,2r)} $. Choose any $y \in g^{-n_0}(x)$ to have a point satisfying the thesis.$\Box$

\begin{lema}\label{come}

   For all $ V $ open set in  $\T^n$ , there exist $y \in V$ and $n_v \in \N$ such that for all $ k \geq n_v: g^k(V) \supset B_{(g^k(y),r)}$.
\end{lema} \pf
For every $V$, Lemma \ref{sucesiondecurvas} gives a point $ y \in V$ and $ n_0 \in \N$ such that $g^k(y)$ is interior to ${g^k(V)}$ and $g^k(y) \in \T^n \setminus B_{(p,2r)}$ for all $k \geq n_0 $. Since $g_{|\T^n \setminus B_{(p,2r)}}$ expands and Theorem \ref{nopasac2} holds in every iterate of $g$ then there exists $ n_V$ such that $ g^k(V) \supset B_{(g^k(y),r)}$ for all $k \geq n_V $. $\Box$

\begin{lema}\label{cubre}
  Given $\varepsilon > 0$ , there exist $ m \in \N$ and a finite family of open balls $\{B_j\}_{j \in \{1,..,d\}}$ in $\T^n$ with $diam(B_j) < \varepsilon $ for all $j$  such that: \begin{itemize}
               \item $\bigcup_{j=1}^d (B_j)= \T^n$,
               \item $ A^m(B_j)= \T^n$ for all $j \in \{1,..,d\}$,
             \item There exists a $C^0$ neighborhood $\mathcal{U_A} $ of $A$ such that for all $ g \in \mathcal{U_A}$ and all $ j \in \{1,..,d\}$ it holds $ g^m(B_j)=\T^n$.
             \end{itemize}
\end{lema}
\pf
The first claim comes from compactness of $\T^n$ and the second from $A$ expands. For the third, define a $C^0$ neighborhood of $A$ explicitly by the formula $\mathcal{U}_A = \{g:\T^n \rightarrow \T^n / \forall p \leq m, \forall j \leq d : |g^p(B_j)-A^p(B_j)| < \varepsilon  \}$ . Endow now every ball in the family with its standard CW complex structure. Since $A$ is linear and $g$ homotopic to $A$, apply Cellular Approximation Theorem to show that $g$ transforms balls in sets that have the homotopy type of balls. Finally, since the boundaries of $g^p(B_j)$ and the boundaries of $A^p(B_j)$ are at distance less or equal than $\ep$ for all $j$, we have that $\mathcal{U}_A$ is the desired neighborhood that completes the proof.$\Box$

\begin{thm}
  $H$ is $C^2$ robustly transitive.
\end{thm}
\pf
Given $\varepsilon = \frac{r}{4}$, choose a neighborhood $\mathcal{U_A}$ in $C^0$ given by Lemma \ref{cubre} and $\mathcal{U}_H \in C^2$ the neighborhood described after Theorem \ref{nopasac2}. Reduce $\mathcal{U}_H$ until it is contained in $\mathcal{U_A}$. Choose an arbitrary $ g \in \mathcal{U}_H$, we prove next that for all $ V$ open set in $\T^n$ there exists $\hat{m} \in \N$ such that $ g^{\hat{m}}(V)= \T^n$.\\ By Lemma \ref{come}, there exist $ y \in V$ and $n_V \in \N$ such that $ g^n(V) \supset B_{(g^n(y),r)}$ for all $n \geq n_V$. Consider now the family of open sets $\{B_j\}_{j \in \{1,..,d\}}$ and $m \in \N$ both given by Lemma \ref{cubre}. Since $diam(B_j)<\frac{r}{4}$ then $g^{n_V}(V) \supset B_j$. The third claim in Lemma \ref{cubre} gives that $g^m(B_j)=\T^n$ which implies that $g^{m+n_V}(V)=\T^n$. This proves that $g$ is transitive. But since $g$ was chosen arbitrarily, the proof holds for all $ g \in \mathcal{U}_H$. Therefore, $H$ is $C^2$ robustly transitive. $\Box$

\section*{Acknowledgements} The author would like to give thanks in the first place to \textrm{Dr. Jorge Iglesias} without whose help this article wouldn't have been possible. And in the second place to \textrm{Dr. Jorge Groisman} and \textrm{Dr. Roberto Markarian} whose thorough reading and insightful comments helped improving the quality of this work.


\begin{thebibliography}{99}

\bibitem[BR]{br}
 P. Berger, A. Rovella.
\newblock On the inverse limit stability of endomorphisms.
\newblock  Ann. Inst. H. Poincar\'e Anal. Non Lin\'eaire , 30.
\newblock p. 463–475.
\newblock 2013.

\bibitem[B]{b}{J.M. Boardman.}{ Singularities of differentiable maps.}
{Publ. Math. IHES 33} {p. 21–57. } {1967.}

\bibitem[BDP]{bdp}
 C. Bonatti, L.J. D\'iaz, E. Pujals.
\newblock A $C^1$ generic dichotomy for diffeomorphisms: weak forms of hyperbolicity or infinitely many sinks or sources.
\newblock  Ann. of Math (2) 158.
\newblock p. 355–418.
\newblock 2003.

\bibitem[GG]{gg}{M. Golubitsky, V. Guillemin.}{ “Stable mappings and their singularities”.}
{Graduate texts in mathematics 14.} {Springer, New York,} {1973.}


\bibitem[KH]{kh} {B. Hasselblatt and A. Katok,}{ “A First Course in Dynamical Systems”.}
{Cambridge University Press, Cambridge,} {2003.}

\bibitem[H]{h}{A. Hatcher.}{“Algebraic Topology”.}
{Cambridge University Press, Cambridge,}{2002.}


\bibitem[ILP]{ilp}
 J. Iglesias, C. Lizana and A. Portela.
\newblock Robust transitivity for endomorphisms admitting critical points.
\newblock  Proc. Amer. Math. Soc. 144, no. 3.
\newblock p. 1235–1250.
\newblock 2016.

\bibitem[IP]{ip}
 J. Iglesias, A. Portela.
\newblock An example of a map which is $C^2$ robustly transitive but not $C^1$ robutstly transitive.
\newblock Colloq. Math. 152, no. 2
\newblock p.285–297.
\newblock 2018.

\bibitem[J]{mj} {M. Jonsson.}{ "Hyperbolic Dynamics of endomorphisms".}
{Royal Institute of Technology, Stockholm,} {1997.}

\bibitem[LP]{lp}
 C. Lizana, E. Pujals.
\newblock Robust transitivity for endomorphisms.
\newblock  Ergod. Th. Dynam. Sys. 33
\newblock p. 1082–1114.
\newblock 2013.

\bibitem[LR1]{lr1}
 C. Lizana, W. Ranter.
\newblock Topological obstructions for robustly transitive endomorphisms on surfaces.
\newblock arXiv:1711.02218
\newblock 2017.

\bibitem[LR2]{lr2}
 C. Lizana, W. Ranter.
\newblock New classes of $C^1$ robustly transitive maps with persistent critical points.
\newblock arXiv:1902.06781
\newblock 2019.

\bibitem[M]{m2}
 R. Ma$\tilde{n}$\'e.
\newblock An ergodic closing lemma.
\newblock  Ann. of Math 116.
\newblock p. 503-540.
\newblock 1982.

\bibitem[P]{p}{F. Przytycki.}{ Anosov Endomorphisms.}
{Studia Mathematica T. LVIII. 3.} {p. 249–285.  }{Polish Academy of Sciences, } {1976.}

\bibitem[S]{s}
M. Sambarino.
\newblock A (short) survey on Dominated Splitting.
\newblock 	arXiv:1403.6050
\newblock 2014.

\bibitem[T]{t}{R. Thom.}{Les singularités des applications différentiables.}
{Ann. Inst. Fourier 6.} {p. 43–87. } {1955.} 

\end{thebibliography}
\end{document}